\documentclass[reqno,12pt]{amsart}
\usepackage{graphicx,amsmath,amssymb,amsthm}

\numberwithin{equation}{section}

\newcommand{\C}{\mathbb C}
\newcommand{\R}{\mathbb R}
\newcommand{\Z}{\mathbb Z}
\newcommand{\N}{\mathbb N}
\renewcommand{\d}{\prime}
\newcommand{\dd}{{\prime \prime}}

\renewcommand{\Re}{{\rm Re}\,}
\renewcommand{\Im}{{\rm Im}\,}
\newtheorem{theorem}{Theorem}[section]
\newtheorem{lemma}[theorem]{Lemma}

\newtheorem{corollary}[theorem]{Corollary}

\newtheorem*{remark}{Remark}

\begin{document}
\title[]
{Trace Formulas for Non-Self-Adjoint Periodic Schr\"odinger Operators \\
and some Applications}
\author[]
{Kwang C.\ Shin}
\address{Department of Mathematics, University of
Missouri, Columbia, MO 65211}
\date{December 10, 2003}
\keywords{Trace formulas, Hill operators, Dirichlet and Neumann eigenvalues}
\footnote{\it 2000 Mathematics Subject Classification. Primary 34A55, 34L15, 34L20.}
\begin{abstract}
Recently, a trace formula for non-self-adjoint periodic Schr\"odinger
operators  in $L^2(\R)$ associated with Dirichlet eigenvalues was proved
 in \cite{FG}. Here we prove a corresponding trace formula associated with
Neumann eigenvalues.

In addition we investigate Dirichlet and Neumann eigenvalues of such
operators. In particular, using the Dirichlet and Neumann trace formulas
we provide detailed information on location of the Dirichlet and Neumann
eigenvalues for the model operator with the potential $Ke^{2ix}$, where $K\in\C$.
\end{abstract}

\maketitle

\baselineskip = 18pt
\section{Introduction}
  Consider the differential expression $L$
\begin{equation}\nonumber
L=-\frac{d^2}{dx^2}+V(x), \quad x\in\R,
\end{equation}
where $V$ is a continuous {\it complex-valued} periodic function on $\R$
of period $\pi$.
The {\it Hill operator} $H$ in $L^2(\R)$ generated by the differential
expression $L$ is defined by
$$(Hf)(x)=L(f(x)),\; x\in\R, \; f\in dom(H)=H^{2,2}(\R), $$
where $H^{2,2}(\R)$ denotes the usual Sobolev space.
Then $H$ is a densely defined closed operator in $L^2(\R)$ (see, e.g.,
\cite[Ch. 5]{MSPE}).

In addition, we define families of differential operators $H_{x_0}^{D}$
and $H_{x_0}^{N}$  in $L^2([x_0,x_0+\pi])$, $x_0\in\R$, as follows:
\begin{align*}
(H_{x_0}^D f)(x)&=L(f(x)),\; x\in [x_0,x_0+\pi], \; f\in dom
(H_{x_0}^D),
\\ (H_{x_0}^N f)(x)&=L(f(x)),\; x\in [x_0,x_0+\pi], \; f\in dom
(H_{x_0}^N),
\end{align*}
where
\begin{align*}
dom (H_{x_0}^D) &= \big\{f\in L^2([x_0,x_0+\pi]) \, \big|\, f,\,f^\d\in
AC([x_0,x_0+\pi]);  \\
& \hspace*{4.2cm}f(x_0+)=0=f(x_0+\pi-)\big\}, \\
dom (H_{x_0}^N) &= \big\{f\in L^2([x_0,x_0+\pi]) \, \big|\, f,\,f^\d\in
AC([x_0,x_0+\pi]);  \\
& \hspace*{4.05cm} f^\d(x_0+)=0=f^\d(x_0+\pi-)\big\}.
\end{align*}

The $L^2(\R)$ spectrum of the Hill operator $H$ is purely continuous and
it is the union of countably many analytic arcs in the complex plane
\cite{FSRB}. On the other hand, the spectra of $H_{x_0}^{D}$  and
$H_{x_0}^{N}$ are purely discrete. We will denote the {\it Dirichlet
eigenvalues} (i.e., the eigenvalues of
$H_{x_0}^{D}$) by $\mu_j(x_0)$, $j\in\N$, and the  {\it Neumann
eigenvalues} (i.e., the eigenvalues of $H_{x_0}^{N}$) by $\nu_k(x_0)$,
$k\in\N_0=\N\cup\{0\}$, where we number these eigenvalues in the
order of non-decreasing magnitudes.

When the potential $V$ is real-valued and periodic, the operators $H$,
$H_{x_0}^D$ and $H_{x_0}^N$ are all self-adjoint, and hence the spectra of
these operators are  subsets of the real line. In this case, the spectrum
of $H$ is a countable union of closed intervals $[E_{2m},\, E_{2m+1}]$,
$m\in\N_0$, on the real line, where $E_m$ are  such that
$L\phi=E_m\phi$ has  a nontrivial (i.e., nonzero) periodic solution of
period $2\pi$. Thus, $E_m$ are eigenvalues of the self-adjoint operator
associated with $L$ under periodic boundary conditions at $0$ and
$2\pi$, and they are all real. Moreover,
\begin{equation}
\nu_0(x_0)\leq E_0, \quad
E_{2m-1}\leq\mu_m(x_0),\,\nu_m(x_0)\leq E_{2m}, \; m\in\N, \; x_0\in\R
\label{int}
\end{equation}
(see, e.g., \cite[Theorem 3.1.1]{MSPE}, \cite{MW}).

Moreover, the following trace formulas hold.
\begin{lemma}\label{real_trace}
Suppose that $V\in C^1(\R)$ is a periodic function on $\R$ and let
$x\in\R$. Then,
\begin{align}
V(x)&=E_0+\sum_{m=1}^{\infty}\big(E_{2m-1}+E_{2m}-2\mu_m(x)\big),
\label{Dt}\\
V(x)&=2\nu_0(x)-E_0+\sum_{m=1}^{\infty}\big(2\nu_m(x)-E_{2m-1}
-E_{2m}\big).\label{Nt}
\end{align}
\end{lemma}
Under the hypothesis that $V\in C^3(\R)$ is a real-valued periodic
function on $\R$, Trubowitz \cite{TRUB} proved the Dirichlet trace formula
(\ref{Dt}). The Neumann trace formula \eqref{Nt} for
real-valued potentials is due to McKean and Trubowitz
\cite{MT}. In 2001, Gesztesy \cite{FG} extended (\ref{Dt}) to
complex-valued periodic $C^1(\R)$ potentials. In Section
\ref{trace_formula}, we will prove the trace formula (\ref{Nt}) for
complex-valued  periodic $C^1(\R)$ potentials, by closely following the
methods in \cite{FG}. We refer to \cite[p.\ 121]{GH03} (cf.\ also
\cite{GS}) for further references and a detailed history of trace
formulas.

In the self-adjoint case, Borg \cite{BORG} showed that $V=E_0$  is the
only real-valued periodic potential, for which the Hill
operator $H$ has the spectrum $[E_0,\infty)$. One can see that if $V\in
C^1(\R)$ is real-valued and  periodic, the Dirichlet trace formula
(\ref{Dt}) yields $V(x)=E_0$ since then
$E_{2m-1}=E_{2m}$ and $E_{2m-1}\leq \mu_m(x)\leq E_{2m}$ for $m\in\N$.
Extensions of this simple observation to reflectionless potentials can be
found in \cite{CGHL}.

Unlike the self-adjoint case, there are infinitely many complex-valued
potentials that generate a half-line spectrum $[0,\infty)$ as
shown by Gasymov \cite{Gasy1} (also, see \cite{Pastur}, \cite{SHIN}). In
these cases, even though all $E_m$ are real and $E_{2m-1}=E_{2m}$
for all $m\in\N$,  the Dirichlet and Neumann eigenvalues $\mu_m(x)$,
$\nu_m(x)$ are not trapped  between $E_{2m-1}$ and $E_{2m}$, due to
non-self-adjointness of the corresponding operators $H_{x}^D$ and
$H_{x}^N$.  From the inverse spectral point of view, this explains the
highly non-uniqueness property of complex-valued periodic potentials
(see, e.g., \cite[p.\ 113]{GH03}). As a concrete example, the potential
$V(x)=e^{2ix}$ generates the spectrum $[0,\,\infty)$ (cf.\ \cite{Gasy1})
and one infers $E_0=0$, $E_{2m-1}=E_{2m}=m^2$, $m\in\N$. Marchenko \cite{March} provides the {\it asymptotic} location of $\mu_j(x)$ and $\nu_k(x)$ (see \eqref{E_asymp}, \eqref{mu} and \eqref{asymp_eq} below). However, up to
date, the precise location of $\mu_j(x)$ and $\nu_k(x)$ for small $j,\,k$, is not known for
this example. The principal application of the trace formulas
proved in this paper will be a determination of the location
of these eigenvalues. This complements the asymptotic results of Marchenko \cite{March}. 

In Section \ref{back_ground}, we will briefly review elements of
Floquet theory and show that algebraic and geometric multiplicities of the
eigenvalues $E_m$ are different for some cases in Theorem
\ref{main_estim}. Our principal new result, the Neumann trace formula for
periodic potentials will be proved in Section \ref{trace_formula}.
In Section \ref{pt_sym}, we will prove the following result:
\begin{theorem}
If $\overline{V(-x)}=V(x)$, $x\in\R$, then for all $x_0\in\R$,
$j\in\N$, and $k\in\N_0$,
$\mu_j(\pi-x_0)=\overline{\mu_j(x_0)}$ and
$\nu_k(\pi-x_0)=\overline{\nu_k(x_0)}$.
\end{theorem}

In Section \ref{deift_work}, we will recall
some unpublished work by Deift \cite{DEIFT} and review
some basic facts on Bessel functions that we use in our final Section
\ref{appl}. In Section \ref{appl}, we will use the trace formulas
(\ref{Dt}) and (\ref{Nt}) to investigate the location
of $\mu_j(x_0)$ and $\nu_k(x_0)$ for the concrete example $V(x)=Ke^{2ix}$,
$K\in\C$. More precisely, we will prove the following theorems.
\begin{theorem}
Suppose that $V(x)=K e^{2ix}$ and $K=|K|e^{2i\varphi_0}$ for some
$\varphi_0\in\R$. Then for all $j\in\N$ and $k\in\N_0$,
$\mu_j(K,x_0)=\mu_j(|K|, x_0+\varphi_0)$ and $\nu_k(K,x_0)=\nu_k(|K|,
x_0+\varphi_0)$.
\end{theorem}

\begin{theorem}\label{main_estim}
Suppose that $V(x)=K e^{2ix}$, $K>0$. Then
\begin{itemize}
\item[(i)] $(j-1)^2\leq\mu_j(0)\leq j^2$ for all $j\in\N$.
\item[(ii)] If $0<K\leq 1$, then  for all $j\in\N$, \\
$j^2-\frac{K}{2}<j^2-\frac{K}{2}+\sum_{m=1}^{j-1}(m^2-\mu_m(0))<\mu_j(0)<j^2$,
and hence $\mu_j(0)\not=E_m$ for all $j\in\N$.
\item[(iii)] If $K>1$ then $M_1=K-2\sum_{m=1}^{\lfloor
\sqrt{K}\rfloor}(m^2-\mu_m(0))>0$ and \\
  $j^2-\frac{M_1}{2}<\mu_j(0)<j^2$ for all $j\geq \sqrt{K}+1$, where
$\lfloor \sqrt{K}\rfloor$ denotes the largest integer that is less than or
equal to $\sqrt{K}$. In particular, $\mu_j(0)\not=E_m$ if $j\geq
\sqrt{K}+1$ and $m\in\N_0$.
\item[(iv)] $k^2\leq\nu_k(0)\leq (k+1)^2$ for all $k\in\N_0$.
\item[(v)] If $0<K\leq 1$, then for all $k\in\N_0$, \\
$k^2<\nu_k(0)<k^2+\frac{K}{2}-\nu_0(0)-\sum_{m=1}^{k-1}(\nu_m(0)
-m^2)<k^2+\frac{K}{2}-\nu_0(0)$,
and hence $\nu_k(0)\not=E_m$ for all $k, m\in\N_0$.
\item[(vi)] If $K>1$ then $M_2=K-2\nu_0(0)-2\sum_{m=1}^{\lfloor
\sqrt{K}\rfloor}(\nu_m(0)-m^2)>0$ and $k^2<\nu_k(0)< k^2+\frac{M_2}{2}$
for all $k\geq \sqrt{K}$. In particular, $\nu_k(0)\not=E_m$ for $k\geq
\sqrt{K}$ and $m\in\N_0$.
\item[(vii)] If $K>0$ then $(2j-1)^2<\mu_{2j-1}(\pi/2)\leq
\mu_{2j}(\pi/2)<(2j)^2$ for all $j\in\N$, and hence $\mu_j(\pi/2)\not=E_m$
for all $j\in\N$ and $m\in\N_0$.
Moreover, if $0<K\leq 1$ then for all $j\in\N$,
$$(2j-1)^2<\mu_{2j-1}(\pi/2)< (2j-1/2)^2<\mu_{2j}(\pi/2)<(2j)^2.$$
\item[(viii)] If $K>0$ then $(2k)^2<\nu_{2k}(\pi/2)\leq
\nu_{2k+1}(\pi/2)<(2k+1)^2$ for all $k\in\N_0$, and hence
$\nu_k(\pi/2)\not=E_m$ for all $k\in\N_0$.
Moreover, if $0<K\leq \frac{1}{2}$ then for all $k\in\N_0$,
$$(2k)^2<\nu_{2k}(\pi/2)< (2k+1/2)^2<\mu_{2k+1}(\pi/2)<(2k+1)^2.$$
\item[(ix)] For every $x_0\in (0,\pi)$, $\mu_j(x_0),\,\nu_k(x_0)\not=E_m$.
\end{itemize}
\end{theorem}

\section{Background}\label{back_ground}
In this section, we introduce some definitions and basic facts on
Floquet theory. In addition, we investigate algebraic and geometric
multiplicities of the eigenvalues $E_m$.

We study
\begin{equation}\label{main_op1}
L\psi(x)=-\psi^\dd(x)+V(x)\psi(x)=\lambda\psi(x),\quad
x\in\R,
\end{equation}
where $\lambda\in\C$ and $V\in L_{\text{loc}}^1(\R)$ is a periodic
function of period $\pi$.

For each $\lambda\in\C$ and $x_0\in\R$, there exists a fundamental
system of solutions $c(\lambda,\,x_0,x),\,
s(\lambda,\,x_0,x)$ of equation (\ref{main_op1}) such that
\begin{align}
c(\lambda,\,x_0,x_0)=1, &
c^\d(\lambda,\,x_0,x_0)=0;\nonumber\\
s(\lambda,\,x_0,x_0)=0, &
s^\d(\lambda,\,x_0,x_0)=1.\nonumber
\end{align}
One can check that these two solutions satisfy the Volterra integral
equations
\begin{align}
c(\lambda,x_0,x)=&\cos[\sqrt{\lambda}(x-x_0)]\nonumber\\
&+\int_{x_0}^xdx_1\frac{\sin[\sqrt{\lambda}(x-x_1)]}{\sqrt{\lambda}}V(x_1)c(\lambda,x_0,x_1),\label{int_eq1}\\
s(\lambda,x_0,x)=&\frac{\sin[\sqrt{\lambda}(x-x_0)]}{\sqrt{\lambda}}\nonumber\\
&+\int_{x_0}^xdx_1\frac{\sin[\sqrt{\lambda}(x-x_1)]}{\sqrt{\lambda}}V(x_1)s(\lambda,x_0,x_1).\label{int_eq2}
\end{align}
{}From these integral equations, along with uniqueness of their
solutions, one can show that
\begin{align}
\frac{\partial}{\partial
x_0}\Big(c(\lambda,x_0,x)\Big)&= (\lambda-V(x_0))s(\lambda,x_0,x),\label{deri_x01}\\
\frac{\partial}{\partial
x_0}\Big(s(\lambda,x_0,x)\Big)&= -c(\lambda,x_0,x).\label{deri_x0}
\end{align}

The {\it monodromy matrix} $M$ associated with \eqref{main_op1} is defined
by
\begin{equation}
M(\lambda,x_0)=\left(
\begin{matrix}
c(\lambda,x_0,x_0+\pi) & s(\lambda,x_0,x_0+\pi)\\
c^\d(\lambda,x_0,x_0+\pi)& s^\d(\lambda,x_0,x_0+\pi)\nonumber
\end{matrix}
\right),
\end{equation}
and the corresponding {\it Floquet discriminant} $\Delta(\lambda)$
(half the trace of $M$) is defined by
\begin{equation}
\Delta(\lambda)=\frac{1}{2}\Big(c(\lambda,\,x_0,x_0+\pi)
+s^\d(\lambda,\,x_0,x_0+\pi)\Big).
\end{equation}
Using (\ref{deri_x01}) and (\ref{deri_x0}), it can be shown that $\Delta$
is independent of $x_0\in\R$.

The Floquet discriminant $\Delta(\lambda)$ is an entire function of order
$\frac{1}{2}$, and $E_m$, $m\in\N_0$, are the zeros of
$\Delta(\lambda)^2-1=0$ (see, e.g., \cite[Ch. 4]{MSPE}). Moreover, since
for each $x_0\in\R$ the entire function
$\lambda\mapsto s(\lambda,x_0,x_0+\pi)$ is of order $\frac{1}{2}$ (see,
e.g., \cite[Ch. 4]{MSPE}, \cite[Ch. 21]{ECT}), the function has
infinitely many zeros $\mu_j(x_0)$, $j\in\N$, by the Picard little theorem
(see, e.g., \cite[Ch.\ 5]{AHLF}).
Similarly, there are infinitely many zeros
$\nu_k(x_0)$, $k\in\N_0$, of the entire
function $c^\d(\lambda,x_0,x_0+\pi)$ of order $\frac{1}{2}$.

\subsection*{Algebraic and geometric multiplicities of the eigenvalues
$E_m$}
  When $V$ is real-valued and periodic, the Hill operator associated with
$L$ in $L^2([0,2\pi])$ with periodic boundary conditions at $0$ and
$2\pi$ (which can easily be defined also for
$V\in L^1_{\text{loc}}(\R)$) is self-adjoint, and hence the algebraic and
geometric multiplicities are the same. However, these two multiplicities
of the eigenvalues $E_m$ do not necessarily agree in the context of
non-self-adjoint Hill operators. In particular, for the concrete example
$V(x)=K e^{2ix}$ that generates a non-self-adjoint Hill operator in
$L^2([0,2\pi])$ with periodic boundary conditions at $0$ and $2\pi$, we
will explain below why these two multiplicities are different for all
$E_m$, $m\in\N$ if $|K|\leq 1$.

The geometric multiplicity of an eigenvalue is the number of linearly
independent eigenfunctions corresponding to the eigenvalue. In addition,
the geometric multiplicity of the eigenvalues $E_m$ corresponding to
periodic or anti-periodic boundary conditions at $0$ and $\pi$  agrees
with the number of linearly independent eigenvectors of the monodromy
matrix $M(E_m,x_0)$ (see, e.g., \cite[Ch.\ 1]{MSPE}).

Since the Wronskian of $c(E_m,x_0,\cdot)$ and $s(E_m,x_0,\cdot)$ is one
(i.e., the determinant of the monodromy matrix equals one), and since the
trace of the monodromy matrix $M(E_m,x_0)$ is $2$ (or $-2$) (i.e.,
$\Delta(E_m)=\pm 1$), we see that $1$ (or $-1$) is the
only eigenvalue of the matrix $M(E_m,x_0)$. In this case, the only way
the monodromy matrix can have two linearly independent eigenvectors occurs
when $M(E_m,x_0)$ equals the $2\times 2$ identity matrix $I_2$ (or
$-I_2$). Thus, the geometric multiplicity of $E_m$ is $2$ if and only if
$c(E_m,x_0,x_0+\pi)=s^\d(E_m,x_0,x_0+\pi)=\pm 1$ and
$c^\d(E_m,x_0,x_0+\pi)=0=s(E_m,x_0,x_0+\pi)$.

In the special case where $V$ is real-valued and periodic, the condition
$E_{2m-1}=E_{2m}$ forces $E_{2m}=\mu_m=\nu_m$, $m\in\N$  (see, e.g.,
\cite[\S 2.3]{MSPE}) and hence
$c^\d(E_{2m},x_0,x_0+\pi)=0=s(E_{2m},x_0,x_0+\pi)$. So
$M(E_{2m}, x_0)=\pm I_2$. Thus, if $E_{2m-1}=E_{2m}$ for some $m\in\N$
then the geometric multiplicities of such $E_{2m}$ are all $2$.
However, when $V$ is complex-valued and periodic, the condition
$E_{2m-1}=E_{2m}$  does not imply
$c^\d(E_{2m},x_0,x_0+\pi)=0=s(E_{2m},x_0,x_0+\pi)$, and hence the
geometric multiplicity of $E_{2m}$ could be 1.

When $V(x)= K e^{2ix}$, $K\in\C$, the Floquet discriminant takes on the
special form $\Delta(\lambda)=\cos(\pi\sqrt{\lambda})$ (see \cite{Gasy1},
\cite{Pastur}, \cite[Theorem 2]{SHIN}). Thus, $E_{2m-1}=E_{2m}=m^2$ for
all $m\in\N$, but as we see from Theorem \ref{main_estim}, if $|K|\leq 1$
then $s(E_{2m},x_0,x_0+\pi)\not=0$ since $\mu_j(x_0)\not=m^2$ for any
$j\in\N$. Thus, if $|K|\leq 1$, the geometric multiplicity of each
$E_m$, $m\in\N$, equals $1$. On the other hand, the algebraic
multiplicity of the eigenvalues $E_m$  could be higher than one. In
general, the algebraic multiplicity of the eigenvalue $E_m$ is the
multiplicity of the zero of
$\Delta(\lambda)\mp 1=0$ at $\lambda=E_m$ (see, e.g., \cite[Theorem
3.3]{FG1}). So except for
$E_0=0$, the eigenvalues $E_{2m-1}=E_{2m}=m^2$ have algebraic
multiplicities $2$ and geometric multiplicities $1$. This result is
mentioned in \cite[p.\ 5]{BB} without proof.

\section{Trace formulas}\label{trace_formula}

In this section, we prove a trace formula for $V$ associated with the
Neumann eigenvalues $\nu_k$, $k\in\N_0$.

In \cite{FG}, Gesztesy proved the following trace formula for Dirichlet
eigenvalues $\mu_j$, $j\in\N$ in the general case where $V\in C^1(\R)$ is
complex-valued and periodic (actually, a larger class of potentials was
considered in \cite{FG}). For simplicity we assume that $V$ has period
$\pi$.

\begin{theorem}\label{trace_diric}
Suppose that the potential $V\in C^1(\R)$ is periodic of period $\pi$.
Then
\begin{equation}\label{Dirichlet_trace1}
V(x)=E_0+\sum_{m=1}^{\infty}\big(E_{2m-1}+E_{2m} -2\mu_m(x)\big), \quad
x\in\R.
\end{equation}
\end{theorem}

Next, we will prove the corresponding trace formula associated with
Neumann eigenvalues $\nu_k(x)$. In Section \ref{appl}, these two trace
formulas will be used to investigate location of  $\mu_j(x)$ and
$\nu_k(x)$ for the example $V(x)= K e^{2ix}$, $K\in\C$.

\begin{theorem}\label{trace_neu}
Suppose that the potential $V\in C^1(\R)$ is periodic of period $\pi$.
Then
\begin{equation}\label{Neumann_trace1}
V(x)=2\nu_0(x)-E_0+\sum_{m=1}^{\infty}\big(2\nu_m(x)
-E_{2m-1}-E_{2m}\big), \quad x\in\R.
\end{equation}
\end{theorem}
\begin{proof}
We will closely follow the proof in the Dirichlet case in \cite{FG}.

Let $G(\lambda,x_1,x_2)$ be the Green's function of the Hill operator $H$
in $L^2(\R)$. Then the diagonal Green's function
$g(\lambda,x)=G(\lambda,x,x)$ of $H$ becomes
\begin{equation}\label{eqeq1}
g(\lambda,x)=-\frac{s(\lambda,x,x+\pi)}{2\sqrt{\Delta^2(\lambda)-1}}
\end{equation}
  (see, e.g., \cite[eq.\ (3.58)]{FG}),
where the branch of the square root is chosen so that the positive real
axis maps onto itself and the square root is analytically continued to the
complex $\lambda$-plane away from the spectrum of $H$.
Next, using (\ref{deri_x0}) and (\ref{eqeq1}), one can show by
straightforward computations that
\begin{equation}\label{eqeq}
\frac{1}{2}g_{xx}(\lambda,x)+(\lambda-V(x))g(\lambda,x) =
\frac{c^\d(\lambda,x,x+\pi)}{2\sqrt{\Delta^2(\lambda)-1}}.
\end{equation}
Together with (\ref{deri_x0}) this implies
\begin{equation}\label{diagonal_asym}
-2g_{xx}(\lambda,x)g(\lambda,x)+g_{x}(\lambda,x)^2+4(V(x)-\lambda)g(\lambda,x)^2
=1.
\end{equation}

Next, we recall that $\{E_m\}_{m\in\N_0}$, $\{\mu_j(x)\}_{j\in\N}$ and
$\{\nu_k(x)\}_{k\in\N_0}$ are the zeros of the functions
$\Delta(\lambda)^2-1$, $s(\lambda,x,x+\pi)$ and $c^\d(\lambda,x,x+\pi)$,
respectively. Thus,
\begin{align}
E_{\substack{2m-1\\
2m}}&\underset{m\to\infty}{=} \left(m+\frac{c_1}{m}+\frac{c_2}{m^2}\pm\frac{\delta_m}{m^2}+\frac{\varepsilon_m^{\pm}}{m^3}\right)^2,\label{E_asymp}\\
\mu_m(x)&\underset{m\to\infty}{=} \left(m+\frac{c_1}{m}
+\frac{c_2}{m^2}+\frac{s_m(x)}{m^2}\right)^2,\label{mu} \\
\nu_m(x)&\underset{m\to\infty}{=} \left(m+\frac{c_1}{m}+\frac{c_2}{m^2}+\frac{c_m(x)}{m^2}\right)^2,
\label{asymp_eq}
\end{align}
where $\{\delta_m\}_{m\in\N}$, $\{\varepsilon^{\pm}_m\}_{m\in\N}$,
$\{s_m(x)\}_{m\in\N}$, $\{c_m(x)\}_{m\in\N} \in \ell^2(\N)$ (see, e.g.,
\cite[Chs. 1, 3]{March}). (Actually, Marchenko \cite{March} did not
explicitly work out the proof of (\ref{asymp_eq}), but using his equation
(1.5.3) along with ideas in the proof of \cite[Theorem 1.5.1]{March},
one can directly prove (\ref{asymp_eq}). When $V$ is real-valued, one can use the
interlacing property \eqref{int} and \eqref{E_asymp}.) 

We write 
\begin{equation}\label{dia_eq2}
g(\lambda,x)=\frac{i\, f_0(\lambda,x)}{2\sqrt{\lambda}}\,\,\text{ for some function $f_0(\lambda, \cdot)\in C^3(\R)$}.
\end{equation}
{}Then  (\ref{E_asymp}) and (\ref{mu}) together with
(\ref{eqeq1}) and the Hardamard factorization theorem  (see, e.g.,
\cite[Ch.\ 5]{AHLF}) imply that
$$f_0(\lambda, x)\underset{\lambda\to i\infty}{=}O(1)\,\,\text{ for each $x\in\R$.}$$ 
{}Next, we infer from
(\ref{diagonal_asym}) and (\ref{dia_eq2}) that $f_0(\lambda, x)^2\underset{\lambda\to i\infty}{=}1$ for every $x\in\R$. Moreover,  one can see that $f_0(\lambda, x)=f_0(x)+f_1(\lambda, x)/\lambda$ for some $f_0(x)^2=1$ and $f_1(\lambda, \cdot)\in C^3(\R)$, where  $f_1(\lambda, x)\underset{\lambda\to i\infty}{=}O(1)$ for each $x\in\R$. 
Thus,
\begin{align}
&g(\lambda,x) = \frac{i\, f_0(x)}{2\sqrt{\lambda}}+\frac{i\, f_1(\lambda,x)}{2(\sqrt{\lambda})^3},\label{dia_eq1}\\
&g_{x}(\lambda,x),\,\, g_{xx}(\lambda,x) \underset{\lambda\to
i\infty}{=} O(|\lambda|^{-3/2}).\label{g_xasy}
\end{align}
Next, from (\ref{diagonal_asym}), (\ref{dia_eq1}) and \eqref{g_xasy} we infer that 
$$f_1(\lambda, x)\underset{\lambda\to i\infty}{=}f_0(x)V(x)/2\,\,\text{ for every $x\in\R$},$$
and hence $g(\lambda,x)$ has the
following asymptotic expression
\begin{equation}\label{dia_eq3}
  g(\lambda,x)\underset{\lambda\to
i\infty}{=} f_0(x)\left(\frac{i}{2\sqrt{\lambda}}+\frac{i}{4(\sqrt{\lambda})^3}V(x)+O(|\lambda|^{-5/2})\right).
\end{equation}
(Actually, it is known that $f_0(x)=1$ but $f_0(x)^2=1$ suffices for
our argument below.) The asymptotic expressions \eqref{g_xasy} and \eqref{dia_eq3} along with (\ref{eqeq})
imply
\begin{equation}\label{asasym}
\left(\frac{c^\d(\lambda,x,x+\pi)}{2\sqrt{\Delta^2(\lambda)-1}}\right)^2
\underset{\lambda\to
i\infty}{=}-\frac{\lambda}{4}+\frac{1}{4}V(x)+O(|\lambda|^{-1}).
\end{equation}

Next, one can choose $n\in\N$  so that the infinite product
\begin{equation}\label{choi_n}
\prod_{m=n}^{\infty}\frac{\nu_m(x)^2}{E_{2m-1}E_{2m}}
\end{equation}
converges absolutely to a non-zero complex number.

If $c^\d(\lambda,x,x+\pi)=0$ for some $x\in\R$ and $\lambda\in\C$, then it
is clear that $t\mapsto c(\lambda,x,t)$ is an eigenfunction of the Neumann
problem in $L^2([x,x+\pi])$ and $\lambda=\nu_k(x)$ for some $k\in\N_0$.
Since $\lambda\mapsto c^\d(\lambda,x,x+\pi)^2$ and
$4(\Delta(\lambda)^2-1)$ are entire functions of order $\frac{1}{2}$,
the Hardamard factorization theorem implies that
\begin{align}
&c^\d(\lambda,x,x+\pi)^2=A_1(x)\prod_{m=0}^{N-1}(\nu_m(x)-\lambda)^2
\prod_{m=N}^{\infty}\left(1-\frac{\lambda}{\nu_m(x)}\right)^2,
\nonumber\\
&4(\Delta(\lambda)^2-1)=B_1\prod_{m=0}^{2N-2}(E_m-\lambda)\prod_{m=N}^{\infty}\left[\left(1-\frac{\lambda}{E_{2m-1}}\right)\left(1-\frac{\lambda}{E_{2m}}\right)\right],\nonumber
\end{align}
for some function $A_1$ on $\R$ and a constant $B_1\in\C$.
Here we chose $N\geq n$  so that  for all $m\geq N$, $E_m\not=0$,
$\nu_m(x)\not=0$ and
\begin{equation}\label{leqeq}
  |\lambda-E_m|\leq |\lambda|\,\text{ for all large $|\lambda|$ with
$\arg(\lambda)=\frac{\pi}{2}$}.
\end{equation}
This is possible because $|\Im(E_m)|\leq \max_{0\leq x\leq\pi}|\Im(V(x))|$
for  $m\in\N_0$.

Thus,
\begin{align}
& \left(\frac{1}{\lambda}\right)\frac{c^\d(\lambda,x,x+\pi)^2}{4(\Delta^2(\lambda)-1)}\nonumber\\
&= -A(x)\frac{(\nu_0(x)-\lambda)^2}{\lambda(E_0-\lambda)}\prod_{m=1}^{\infty}\left[\frac{(\nu_m(x)-\lambda)^2}{(E_{2m-1}-\lambda)(E_{2m}-\lambda)}\right]\nonumber\\
&= A(x)\frac{\left(1-\frac{\nu_0(x)}{\lambda}\right)^2}{1-\frac{E_0}{\lambda}}\prod_{m=1}^{\infty}\left[\frac{\left(1-\frac{\nu_m(x)}{\lambda}\right)^2}{\left(1-\frac{E_{2m-1}}{\lambda}\right)\left(1-\frac{E_{2m}}{\lambda}\right)}\right],\label{asymp_eq1}
\end{align}
where
\begin{equation}
A(x)=-\frac{A_1(x)}{B_1}\prod_{m=N}^{\infty}\frac{E_{2m-1}E_{2m}}{\nu_m(x)^2}.
\end{equation}

Next, we have
\begin{align}
& \ln\left(\prod_{m=1}^{N}\left[\frac{(\nu_m(x)-\lambda)^2}{(E_{2m-1}-\lambda)(E_{2m}-\lambda)}\right]\right)\nonumber\\
&= \sum_{m=1}^{N}\left[2\ln\left(1-\frac{\nu_m(x)}{\lambda}\right)-\ln\left(1-\frac{E_{2m-1}}{\lambda}\right)-\ln\left(1-\frac{E_{2m}}{\lambda}\right)\right]\nonumber\\
&\underset{\lambda\to
i\infty}{=} -\frac{1}{\lambda}\sum_{m=1}^{N}\left[2\nu_m(x)-E_{2m-1}-E_{2m}\right]+O(|\lambda|^{-2}),\label{first_nterm}
\end{align}
where we used $\ln(1-t)\underset{t\to 0}=-t+O(t^2).$
Moreover,
\begin{align}
&\sum_{m=N+1}^{\infty}\left|\frac{(2\nu_m(x)-E_{2m-1}-E_{2m})}{\left(\frac{E_{2m-1}}{\lambda}-1\right)\left(\frac{E_{2m}}{\lambda}-1\right)}\right|\nonumber\\
&\leq\sum_{m=N+1}^{\infty}\frac{\left|2\nu_m(x)-E_{2m-1}-E_{2m}\right|}{|\lambda|^{-2}\left|\left(E_{2m-1}-\lambda\right)\left(E_{2m}-\lambda\right)\right|}
\nonumber\\
&\leq\sum_{m=N+1}^{\infty}\frac{\left|2\nu_m(x)-E_{2m-1}-E_{2m}\right|}{|\lambda|^{-2}\left|\lambda^2\right|}<+\infty,
\nonumber
\end{align}
and  by (\ref{E_asymp}), (\ref{asymp_eq}) and (\ref{leqeq}),
\begin{align}
&\sum_{m=N+1}^{\infty}\frac{\left|\nu_m(x)^2
-E_{2m-1}E_{2m}\right|}{\left|\left(\frac{E_{2m-1}}{\lambda}-1\right)
\left(\frac{E_{2m}}{\lambda}-1\right)\right|} \nonumber \\
&\leq\sum_{m=N+1}^{\infty}\left|\nu_m(x)^2-E_{2m-1}E_{2m}\right|<+\infty.
\nonumber
\end{align}
Thus, we obtain
\begin{align}
&\ln\left(\prod_{m=N+1}^{\infty}\left[\frac{(\nu_m(x)-\lambda)^2}{(E_{2m-1}-\lambda)(E_{2m}-\lambda)}\right]\right)\nonumber\\
&=\sum_{m=N+1}^{\infty}\ln\left[1-\frac{\left(E_{2m-1}-\lambda\right)\left(E_{2m}-\lambda\right)-(\nu_m(x)-\lambda)^2}{\left(E_{2m-1}-\lambda\right)\left(E_{2m}-\lambda\right)}\right]\nonumber\\
&=\sum_{m=N+1}^{\infty}\ln\left[1-\frac{(2\nu_m(x)-E_{2m-1}-E_{2m})\lambda+(E_{2m-1}E_{2m}-\nu_m(x)^2)}{\lambda^2\left(\frac{E_{2m-1}}{\lambda}-1\right)\left(\frac{E_{2m}}{\lambda}-1\right)}\right]\nonumber\\
&=-\sum_{m=N+1}^{\infty}\left[\frac{(2\nu_m(x)-E_{2m-1}-E_{2m})\lambda+(E_{2m-1}E_{2m}-\nu_m(x)^2)}{\lambda^2\left(\frac{E_{2m-1}}{\lambda}-1\right)\left(\frac{E_{2m}}{\lambda}-1\right)}\right]\nonumber\\
&+O\left(\sum_{m=N+1}^{\infty}\left|\frac{(2\nu_m(x)-E_{2m-1}-E_{2m})\lambda+(E_{2m-1}E_{2m}-\nu_m(x)^2)}{\lambda^2\left(\frac{E_{2m-1}}{\lambda}-1\right)\left(\frac{E_{2m}}{\lambda}-1\right)}\right|^2\right)\nonumber\\
&=-\frac{1}{\lambda}\sum_{m=N+1}^{\infty}\frac{(2\nu_m(x)-E_{2m-1}-E_{2m})}{\left(\frac{E_{2m-1}}{\lambda}-1\right)\left(\frac{E_{2m}}{\lambda}-1\right)}\nonumber\\
&-\frac{1}{\lambda^2}\sum_{m=N+1}^{\infty}\frac{(E_{2m-1}E_{2m}-\nu_m(x)^2)}{\left(\frac{E_{2m-1}}{\lambda}-1\right)\left(\frac{E_{2m}}{\lambda}-1\right)}+O(|\lambda|^{-2}),\label{estim_eq}
\end{align}
where we used again (\ref{leqeq}) to bound the error term. Next, we use
the dominated convergence theorem (with the discrete counting measure) in
(\ref{estim_eq}) to derive
\begin{align}
&\ln\left(\prod_{m=N+1}^{\infty}\left[\frac{(\nu_m(x)
-\lambda)^2}{(E_{2m-1}-\lambda)(E_{2m}-\lambda)}\right]\right)\nonumber\\
&\underset{\lambda\to
i\infty}{=}-\frac{1}{\lambda}\sum_{m=N+1}^{\infty}(2\nu_m(x)
-E_{2m-1}-E_{2m})+O\left(\frac{1}{|\lambda|^2}\right). \label{final_eq}
\end{align}
Finally,  (\ref{asymp_eq1}), (\ref{first_nterm}) and (\ref{final_eq})
along with the fact $e^{t}\underset{t\to 0}{=}1+O(t)$ yield
\begin{align}
&\frac{1}{\lambda}\frac{c^\d(\lambda,x,x+\pi)^2}{4(\Delta^2(\lambda)-1)}
\underset{\lambda\to
i\infty}{=} A(x)\Big[1-\frac{2\nu_0(x)-E_0}{\lambda}\nonumber\\
& \quad -\frac{1}{\lambda}\sum_{m=1}^{\infty}
\left[2\nu_m(x)-E_{2m-1}-E_{2m}\right]
+O\left(\frac{1}{|\lambda|^2}\right)\Big]. \label{t}
\end{align}
Thus, comparing \eqref{t} with (\ref{asasym}), we infer
$A(x)=-\frac{1}{4}$ and the desired trace formula \eqref{Neumann_trace1}.
\end{proof}

\section{The case when $\overline{V(-x)}=V(x)$}\label{pt_sym}

In this section, we prove some general results on the location of
Dirichlet and Neumann eigenvalues.

The following lemma on periodicity of the eigenvalues is well-known and
hence we omit the proof.

\begin{lemma}\label{lemma_peri}
Suppose that $V\in L^1_{\text{loc}}(\R)$ is periodic of period $\pi$. Then
$\mu_j(x_0)$, $j\in\N$, and $\nu_k(x_0)$, $k\in\N_0$, are all periodic
functions of period $\pi$.
\end{lemma}

Next, we prove the following theorem, regarding a certain symmetry of the
eigenvalues $\mu_j(x_0)$ and $\nu_k(x_0)$.
\begin{theorem}
Suppose that $V\in L^1_{\text{loc}}(\R)$ is periodic of period $\pi$
and that $\overline{V(-x)}=V(x)$ for all $x\in\R$, or equivalently, $\overline{V(\pi-x)}=V(x)$ for all $x\in[0,\pi]$. 
Then for every  $j\in\N$, $k\in\N_0$, and $x_0\in[0,\pi]$,
\begin{equation}
\overline{\mu_j(\pi-x_0)}=\mu_j(x_0)\quad\text{and}\quad
\overline{\nu_k(\pi-x_0)}=\nu_k(x_0).
\end{equation}
\end{theorem}
\begin{proof}
Let $y(\mu_j(x_0),x)=s(\mu_j(x_0), x_0,x)$.  Then since
$\overline{V(-x)}=V(x)$, we see that $\overline{y(\mu_j(x_0),-x)}$ is also
a solution of
\begin{equation}\label{main_eq12}
-\psi^\dd(x)+V(x)\psi(x)=\mu_j(x_0)\psi(x)
\end{equation}
with $\mu_j(x_0)$ replaced by $\overline{\mu_j(x_0)}$. (To see this, we
take the complex conjugate of equation (\ref{main_eq12}) and replace
$x$ by $-x$.) We write
$\overline{y(\mu_j(x_0),-x)}=y(\overline{\mu_j(x_0)},x)$.

Thus, $y(\mu_j(x_0),x_0)=0=y(\mu_j(x_0),x_0+\pi)$ yields
$$y(\overline{\mu_j(x_0)},-x_0-\pi)=0=y(\overline{\mu_j(x_0)},-x_0)=y(\overline{\mu_j(x_0)},(-x_0-\pi)+\pi).$$
This implies
$$\mu_j(-x_0-\pi)=\overline{\mu_j(x_0)}.$$
By Lemma \ref{lemma_peri} we then have
$$\mu_j(\pi-x_0)=\mu_j((\pi-x_0)-2\pi)=\overline{\mu_j(x_0)}.$$

Similarly, we can show that
$\nu_k(\pi-x_0)=\overline{\nu_k(x_0)}$.
\end{proof}

Thus, for all $x_0\in[0,\pi/2]$, $\mu_j(x_0)$ and $\nu_k(x_0)$ contain all
information of the eigenvalues for $x_0\in\R$. Moreover, we have the
following corollary regarding reality of the eigenvalues $\mu_j(x_0)$,
$\nu_k(x_0)$ for $x_0=0$, $\frac{\pi}{2}$.
\begin{corollary}\label{cor_real}
Suppose that $V\in L^1_{\text{loc}}(\R)$ is periodic of period $\pi$
and that $\overline{V(-x)}=V(x)$ for all $x\in\R$. Then for all $j\in\N$,
$k\in\N_0$, $\mu_j(0)$, $\mu_j(\frac{\pi}{2})$, $\nu_k(0)$ and
$\nu_k(\frac{\pi}{2})$ are all real.
\end{corollary}
\begin{proof}
Since
$\mu_j(0)=\mu_j(\pi)=\overline{\mu_j(\pi-\pi)}=\overline{\mu_j(0)}$,
$\mu_j(0)$ are all real. Similarly, $\nu_k(0)$ are all real.
Moreover,
$$\mu_j(\frac{\pi}{2})=\overline{\mu_j(\pi-\frac{\pi}{2})}
=\overline{\mu_j(\frac{\pi}{2})}, \quad
\nu_k(\frac{\pi}{2})=\overline{\nu_k(\pi-\frac{\pi}{2})}
=\overline{\nu_k(\frac{\pi}{2})}.$$
\end{proof}

\section{Bessel functions and some unpublished work of
Deift}\label{deift_work}
In this section, we reconsider some results in
an unpublished manuscript by Deift \cite{DEIFT}, where he
explicitly expressed $c^\d(\lambda,x_0,x_0+\pi)$ and
$s(\lambda,x_0,x_0+\pi)$ for $V(x)=K e^{2ix}$, $K\in\C$, in terms of
Bessel functions. In addition we will introduce some useful facts on
Bessel functions.

Consider the Schr\"odinger equation
\begin{equation}\label{main_eq}
-\psi^\dd(x)+K e^{2ix}\psi(x)=\lambda\psi(x),\quad
x\in\R,
\end{equation}
where $\lambda,\, K\in\C$.
  Gasymov \cite{Gasy1} (also, see \cite{Pastur}, \cite[Theorem 2]{SHIN})
showed that if $V(x)=K e^{2ix}$, $K\in\C$, then
$$\Delta(\lambda)=\cos(\pi\sqrt{\lambda}).$$
{}From this fact it is clear that
\begin{equation}
  E_0=0,\quad E_{2m-1}=E_{2m}=m^2,\quad m\in\N.
\end{equation}

Next, the Bessel function $J_{\nu}(u)$ of the first kind
is given by
\begin{equation}\label{J_def}
J_{\nu}(u)=\left(\frac{u}{2}\right)^{\nu}
\sum_{m=0}^{\infty}\frac{(-1)^m}{m!\Gamma(m+\nu+1)}
\left(\frac{u}{2}\right)^{2m},\quad
u, \nu \in\C,
\end{equation}
where we choose the negative real axis as its branch cut for
$\nu\not\in\Z$. The Bessel function $J_{\nu}$ solves the following
differential equation
\begin{equation}\label{bessel_eq}
\frac{d^2}{du^2}J_{\nu}(u)+\frac{1}{u}\frac{d}{du}J_{\nu}(u)+\frac{u^2-\nu^2}{u^2}J_{\nu}(u)=0.
\end{equation}
One can show that
$x\mapsto J_{\sqrt{\lambda}}(\sqrt{K}e^{ix})$ is a solution of
(\ref{main_eq}) (see, e.g., \cite[p.\ 196]{Fl}). Moreover,
\begin{align*}
c(\lambda,0,x)&=\frac{\pi\sqrt{K}}{2}
\left(Y_{\sqrt{\lambda}}^\d(\sqrt{K})J_{\sqrt{\lambda}}
(\sqrt{K}e^{ix})-
J_{\sqrt{\lambda}}^\d(\sqrt{K})
Y_{\sqrt{\lambda}}(\sqrt{K}e^{ix})\right), \\
s(\lambda,0,x)&=\frac{\pi
i}{2}\left(Y_{\sqrt{\lambda}}(\sqrt{K})
J_{\sqrt{\lambda}}(\sqrt{K}e^{ix})-
J_{\sqrt{\lambda}}(\sqrt{K})Y_{\sqrt{\lambda}}
(\sqrt{K}e^{ix})\right),
\end{align*}
where $Y_{\nu}$ is the Bessel function of the second kind (also a
solution of \eqref{bessel_eq}), defined by
\begin{equation}\label{second_bessel}
Y_{\nu}(u)=\frac{J_{\nu}(u)\cos(\nu\pi)-J_{-\nu}(u)}{\sin(\nu\pi)}.
\end{equation}

In his unpublished manuscript \cite{DEIFT}, Deift showed  that
\begin{align}
s(\lambda,0,\pi)&= \pi
J_{\sqrt{\lambda}}(\sqrt{K})J_{-\sqrt{\lambda}}(\sqrt{K}),\label{s_eq}\\
c^\d(\lambda, 0, \pi)&= \pi K
J_{\sqrt{\lambda}}^\d(\sqrt{K})J_{-\sqrt{\lambda}}^\d(\sqrt{K}),\label{c_eq}
\end{align}
where he used  (\ref{second_bessel}) with
\begin{align}
J_{\nu}(ue^{i\pi})&= e^{\nu\pi i}J_{\nu}(u),\nonumber\\
Y_{\nu}(ue^{i\pi})&= e^{-\nu\pi
i}Y_{\nu}(u)+2i\cos(\nu\pi)J_{\nu}(u).\nonumber
\end{align}
In Section \ref{appl}, we will extensively use equations (\ref{s_eq}) and
(\ref{c_eq}), along with the trace formulas \eqref{Dirichlet_trace1},
\eqref{Neumann_trace1} to investigate the location of
$\mu_j$ and $\nu_k$.

Next, we list a number of basic facts on Bessel functions
that will be used in the next section.
\begin{lemma}\label{lemma_bessel} ${}$\\
\vspace*{-5mm}
\begin{itemize}
\item[(1)] If $\nu$ is real, then $J_{\nu}(u)$ is real
for all $u\in\R$.
\item[(2)] If $\nu\geq -1$, then all zeros of $J_{\nu}(u)$ are real and
if
$\nu\geq 0$, then all zeros of $J_{\nu}^\d(u)$ are real.
\item[(3)] If $\nu\geq 0$, then the smallest positive
zeros of $J_{\nu}(u)$ and $J_{\nu}^\d(u)$ are greater
than $\nu$.
\item[(4)] If $n\in\Z$, then
$J_{-n}(u)=(-1)^{n}J_n(u)$  and
$J_{-n}^\d(u)=(-1)^{n}J_n^\d(u)$, $\,u\in\C$.
\item[(5)] $J_{\nu}$ and $J_{-\nu}$ are linearly independent if and only
if $\nu\not\in\Z$.
\item[(6)]  $J_{n}$ and $J_{m}$ do not have a
common zero if $m,n\in\N_0$ with $m\not=n$. 
\item[(7)] $J_{n-1}^\d$ and $J_{n}^\d$ do not have a common zero for
$n\in\N$.
\end{itemize}
\end{lemma}
\begin{proof}
See, for example, \cite[Ch.\ 15]{WATSON} and \cite[Ch.\ 9]{MAIS} for
proofs of these results. We also note that assertions (1), (4) and (5)
can be derived directly from (\ref{J_def}).
\end{proof}

\section{Applications of the trace formulas}\label{appl}
In this section, we will prove a number of results regarding the location
of the Dirichlet and Neumann eigenvalues for the potential $V(x)=K
e^{2ix}$, the collection of which becomes Theorem \ref{main_estim}.

In the following theorem, we will use $\mu_j(K,x_0)$ and
$\nu_k(K,x_0)$ to explicitly indicate the $K$-dependence of these
eigenvalues.
\begin{theorem}\label{thm_trans}
Suppose that $V(x)=K e^{2ix}$ and $K=|K|e^{2i\varphi_0}$ for some
$\varphi_0\in\R$. Then for each $x\in \R$, $j\in\N$, and $k\in\N_0$,
$$\mu_j(K,x)=\mu_j(|K|, x+\varphi_0) \text{ and }
\nu_k(K,x)=\nu_k(|K|, x+\varphi_0).$$
\end{theorem}
\begin{proof}
With $K=|K|e^{2i\varphi_0}$, (\ref{main_eq}) becomes
\begin{equation}\label{main_eq1}
-\psi^\dd(x)+|K|e^{2i(x+\varphi_0)}\psi(x)=\lambda\psi(x),\quad
x\in\R.
\end{equation}
Next, we consider the equation
\begin{equation}\label{main_eq2}
-\psi^\dd(x)+|K|e^{2ix}\psi(x)=\lambda\psi(x),\quad
x\in\R.
\end{equation}
Clearly, $y(x)$ is a solution of (\ref{main_eq1}) if
and only if $y_1(x)=y(x-\varphi_0)$ is a solution of
(\ref{main_eq2}). Moreover, $y(x)=0=y(x+\pi)$ if and only if
$y_1(x+\varphi_0)=0=y_1(x+\varphi_0+\pi)$. This
implies that
$\mu_j(|K|e^{2i\varphi_0},x)=\mu_j(|K|,x+\varphi_0)$, $j\in\N$.
Similarly one proves that
$\nu_k(|K|e^{2i\varphi_0},x)=\nu_k(|K|,x+\varphi_0)$, $k\in\N_0$.
\end{proof}
\begin{remark}
{\rm (i) From the proof of Theorem \ref{thm_trans}, we conclude that for
all $x\in\R$, $j\in\N$, and $k\in\N_0$,
\begin{equation}\label{trans_eq}
\mu_j(K,x)=\mu_j(e^{2ix}K,0) \text{ and }
\nu_k(K,x)=\nu_k(e^{2ix}K,0).
\end{equation}

(ii) By Theorem \ref{thm_trans}, the periodic curves in the complex plane,
generated by $x\mapsto\mu_j(x)$ and $x\mapsto\nu_k(x)$ remain the
same for all $K\in\C$ with the same magnitude. Thus, we will focus on the
case $K>0$.}
\end{remark}
Next, we will provide more precise location of the Dirichlet eigenvalues
$\mu_j(0)$.
\begin{theorem}\label{mu_real}
Suppose that $V(x)=K e^{2ix}$ and $K>0$. Then
\begin{itemize}
\item[(i)] $(j-1)^2\leq\mu_j(0)\leq j^2$ for all $j\in\N$.
\item[(ii)] If $0<K\leq 1$, then  for all $j\in\N$, \\
$j^2-\frac{K}{2}<j^2-\frac{K}{2}+\sum_{m=1}^{j-1}(m^2-\mu_m(0))<\mu_j(0)<j^2$,
and hence $\mu_j(0)\not=E_m$ for all $j\in\N$, $m\in\N_0$.
\item[(iii)] If $K>1$ then $M_1=K-2\sum_{m=1}^{\lfloor
\sqrt{K}\rfloor}(m^2-\mu_m(0))>0$ and \\
  $j^2-\frac{M_1}{2}<\mu_j(0)<j^2$ for all $j\geq \sqrt{K}+1$, where
$\lfloor \sqrt{K}\rfloor$ denotes the largest integer that is less than or
equal to $\sqrt{K}$. In particular, $\mu_j(0)\not=E_m$ if $j\geq
\sqrt{K}+1$ and $m\in\N_0$.
\end{itemize}
\end{theorem}
\begin{proof}[Proof of (ii)]
We recall that $\mu_j(0)$, $j\in\N$, are the zeros of $s(\lambda,0,\pi)$
and from (\ref{s_eq}),
\begin{equation}\label{J_eq1}
s(\lambda,0,\pi)=\pi
J_{\sqrt{\lambda}}(\sqrt{K})J_{-\sqrt{\lambda}}(\sqrt{K}).
\end{equation}
Also, by Corollary \ref{cor_real}, $\mu_j(0)\in\R$ for all $j\in\N$.

Below, we will show the existence of $\mu_j(0)$ in certain intervals
applying the intermediate value theorem. Subsequently the convergence of
the series in the trace formula in Theorem \ref{trace_diric} will be used
to show that there exist no other $\mu_j(0)$.

First, we will show that if $0<K\leq 1$, then $s((2n)^2,0,\pi)>0$ and
$s((2n-1)^2,0,\pi)<0$, $n\in\N$. The continuous function
$s(\lambda,0,\pi)$ is real for
$\lambda\geq 0$ (in fact, it is real for all
$\lambda\in\R$ because if $\lambda<0$ then the two
Bessel functions in (\ref{J_eq1}) at $\sqrt{K}\in\R$ are complex
conjugates of each other since $\Re(\sqrt{\lambda})=0$).

By Lemma \ref{lemma_bessel} (4),
$$J_{-n}(u)=(-1)^nJ_{n}(u), \quad  n\in\Z$$
and hence
\begin{equation}
s(n^2,0,\pi)=(-1)^n\left(J_n(\sqrt{K})\right)^2.
\end{equation}
By Lemma \ref{lemma_bessel} (2) and (3), for each $\sqrt{\lambda}\geq
0$, the zeros of $u\mapsto J_{\sqrt{\lambda}}(u)$ are all real, and
positive zeros of these functions are all greater than
$\sqrt{\lambda}$. Thus, if $0<K\leq 1$, then for every
$n\geq 1$, $J_n(\sqrt{K})\not=0$, and
$J_0(\sqrt{K})\not=0$ by (\ref{J_def}) since the terms in the series
(\ref{J_def}) have alternating signs and since absolute values of
these terms are strictly decreasing if $\sqrt{K}\leq 2$. Thus, the
sequence $n \mapsto s(n^2,0,\pi)=(-1)^n\left(J_n(\sqrt{K})\right)^2$
alternates its sign for all $n\geq 0$. By the intermediate value theorem
at least one of $\mu_j(0)$ lies in every open interval $((n-1)^2,n^2)$,
$n\in\N$. Next, using the trace formula \eqref{Dirichlet_trace1}, we
will show that there is precisely one $\mu_j(0)$ in every interval
$((n-1)^2,n^2)$.

The trace formula (\ref{Dirichlet_trace1}) at $x_0=0$ reads
\begin{align}
K=Ke^{2i
x_0}\Big|_{x_0=0}&= E_0+\sum_{j=1}^{\infty}\left(E_{2j-1}+E_{2j}-2\mu_j(0)\right)\nonumber\\
&= 0+\sum_{j=1}^{\infty}\left(2j^2-2\mu_j(0)\right).\label{trace_for}
\end{align}
This implies that if there were more than one $\mu_j(0)$ in
some $((n-1)^2,n^2)$ or if there were one $\mu_j(0)$ on the negative real axis, then $j^2-\mu_j(0)>1$ for all $j\geq n+1$, and hence 
the sum in (\ref{trace_for}) would be divergent. 
Thus, there is exactly one $\mu_j(0)$ in every interval $((n-1)^2,n^2)$,
$n\in\N$. Since the sum is $K$ and since $2j^2-2\mu_j(0)>0$ for all
$j\in\N$, we have
$$j^2-\frac{K}{2}<j^2-\frac{K}{2}+\sum_{m=1}^{j-1}\mu_m(0)
<\mu_j(0)<j^2, \quad j\in\N.$$

\noindent{\it Proof of (i) and (iii).}
There is at least one $\mu_j(0)$ in each closed interval
$[(n-1)^2,n^2]$, $n\in\N$. Otherwise, $s((n-1)^2,0,\pi)\not=0$ and
$s(n^2,0,\pi)\not=0$. By (\ref{J_eq1}), there would be at least
one $\mu_j(0)$ in the interior of the interval by the intermediate value
theorem. Next, since the smallest positive zero of $J_{\nu}(u)$, $\nu\geq
0$, is greater than $\nu$, the sequence $n\mapsto  s(n^2,0,\pi)$
alternates its sign for  $n\geq\sqrt{K}$. However, it is possible that $
s(n^2,0,\pi)=0$  for some $0\leq n<\sqrt{K}$.

Suppose that $s(n_0^2,0,\pi)=0$ for some $n_0\in\N$, $0<n_0<\sqrt{K}$.
Then by (\ref{J_eq1}) either $J_{n_0}(\sqrt{K})=0$ or
$J_{-n_0}(\sqrt{K})=0$. However, since  $J_{n_0}$ and $J_{-n_0}$ are
linearly dependent, by Lemma \ref{lemma_bessel} (4), one concludes that
$J_{n_0}(\sqrt{K})=0$ and $J_{-n_0}(\sqrt{K})=0$. Next, we will show that
$s(\lambda,0,\pi)$ has at least double zeros at $\lambda=n_0^2$ if
$s(n_0^2, 0, \pi)=0$.

Since $J_{n_0}(\sqrt{K})=0$, one infers that  $\sqrt{\lambda}\mapsto
J_{\sqrt{\lambda}}(\sqrt{K})=(\sqrt{\lambda}-n_0)f(\sqrt{\lambda})$ for
some entire function $f$. Since $J_{-n_0}(\sqrt{K})=0$, also
$J_{\sqrt{\lambda}}(\sqrt{K})
=(\sqrt{\lambda}-n_0)(\sqrt{\lambda}+n_0)f_1(\lambda)$ for some entire
function $f_1$. Thus,
\begin{align}
s(\lambda,0,\pi)&= \pi
J_{\sqrt{\lambda}}(\sqrt{K})J_{-\sqrt{\lambda}}(\sqrt{K})\nonumber\\
&= \pi (\lambda-n_0^2)^2 f_1(\sqrt{\lambda})f_1(-\sqrt{\lambda}) \nonumber
\end{align}
and hence $s(\lambda,0,\pi)$ has at least a double zero at
$\lambda=n_0^2$. If $n_0=0$ then since $n_0=-n_0$, we have for
some entire function $f_1$,
$$s(\lambda,0,\pi)=\pi \lambda f_1(\sqrt{\lambda})f_1(-\sqrt{\lambda}),$$
and hence the zero can be simple or it can be of higher order.

By Lemma \ref{lemma_bessel} (6), if  $s(n_0^2,0,\pi)=0$ for some integer
$0< n_0<\sqrt{K}$, then $s(n^2,0,\pi)\not=0$ for $n\in\N_0$ with
$n\not=n_0$.
If $n_0=0$, then the zero must be simple, due to the convergence of the
sum in  (\ref{trace_for}) since $(j-1)^2<\mu_j(0)<j^2$ for all $j\geq 2$.
If $n_0\not=0$ then the zero is not simple from the above argument. In
fact, the zero must be of order $2$; otherwise, the sum in
(\ref{trace_for}) would diverge.

If $s(n^2,0,\pi)\not=0$ for all $n\in\N_0$, then  $(j-1)^2<\mu_j(0)<j^2$
for all $j\in\N$ like in the case (ii). Thus, we proved (i).

The proof of (iii) is analogous to that of (ii).
\end{proof}

Next, we study the location of Neumann eigenvalues $\nu_k(0)$.
\begin{theorem}
Suppose that $V(x)=K e^{2ix}$ and $K>0$. Then
\begin{itemize}
\item[(i)] $k^2\leq\nu_k(0)\leq (k+1)^2$ for all $k\in\N_0$.
\item[(ii)] If $0<K\leq 1$, then for all $k\in\N_0$, \\
$k^2<\nu_k(0)<k^2+\frac{K}{2}-\nu_0(0)
-\sum_{m=1}^{k-1}(\nu_m(0)-m^2)<k^2+\frac{K}{2}-\nu_0(0)$,
and hence $\nu_k(0)\not=E_m$ for $k\in\N_0$, $m\in\N_0$.
\item[(iii)] If $K>1$ then $M_2=K-2\nu_0(0)-2\sum_{m=1}^{\lfloor
\sqrt{K}\rfloor}(\nu_m(0)-m^2)>0$ and $k^2<\nu_k(0)< k^2+\frac{M_2}{2}$
for all $k\geq \sqrt{K}$. In particular, $\nu_k(0)\not=E_m$ for $k\geq
\sqrt{K}$, $m\in\N_0$.
\end{itemize}
\end{theorem}
\begin{proof}[Proof of (ii)]
The arguments are very similar to that in the proof of Theorem
\ref{mu_real} (ii). However, in the trace formula (\ref{Neumann_trace1}),
$\nu_0(x_0)$ is paired with $E_0$, unlike in the Dirichlet case
(\ref{Dirichlet_trace1}). Hence, a more careful analysis is needed in the
present Neumann case.

We recall that $\nu_k(0)$, $k\in\N_0$, are the zeros of
\begin{equation}
c^\d(\lambda,0,\pi)=\pi K
J_{\sqrt{\lambda}}^\d(\sqrt{K})J_{-\sqrt{\lambda}}^\d(\sqrt{K}),
\end{equation}
a real-valued continuous function on the real line. Moreover, by Corollary
\ref{cor_real}, $\nu_k(0)\in\R$ for all $k\in\N_0$.

If $0<K\leq 1$, as in the case of Dirichlet eigenvalues, the sequence
$n \mapsto
c^\d(n^2,0,\pi)=(-1)^n\left(J_n^\d(\sqrt{K})\right)^2$
alternates its sign for all $n\geq 1$. By (\ref{J_def}), we see that
$J_0^\d(\sqrt{K})\not=0$ as in the case of $J_0(\sqrt{K})\not=0$ for
$\sqrt{K}\leq 2$.
Hence, by the intermediate value theorem there is at least one
$\nu_k(0)$ in  the  interval $((n-1)^2,n^2)$ for every $n\in\N$. Moreover,
$\nu_k(0)\not=n^2$ for any $n\in\Z$. Next, we will show that
$\nu_0(0)>0$, that is, there is no negative Neumann eigenvalue since
$c^\d(0,0,\pi)\not=0$.

The trace formula (\ref{Neumann_trace1}) at $x_0=0$ reads
\begin{align}
K=Ke^{2i
x_0}\Big|_{x_0=0}&= 2\nu_0(0)-E_0+
\sum_{k=1}^{\infty}\left(2\nu_k(0)-E_{2k-1}-E_{2k}\right)\nonumber\\
&= 2\nu_0(0)+\sum_{k=1}^{\infty}\left(2\nu_k(0)-2k^2\right).
\label{trace_for1}
\end{align}

Next, suppose that $J_{\sqrt{\lambda}}^\d(\sqrt{K})=0$ at
$\lambda=\nu_0(0)<0$. Then since $J_{\sqrt{\lambda}}^\d(\sqrt{K})$ is an
entire function of $\sqrt{\lambda}$, we can write
$$J_{\sqrt{\lambda}}^\d(\sqrt{K})=(\sqrt{\lambda}-i\sqrt{|\nu_0(0)|})f(\sqrt{\lambda})$$
  for some entire function $f$.
  Since
$$J_{-\sqrt{\nu_0(0)}}^\d(\sqrt{K})=\overline{J_{\sqrt{\nu_0(0)}}^\d(\sqrt{K})}
=0,$$
  we see that $f(-i\sqrt{|\nu_0(0)|})=0$, and hence
$$J_{\sqrt{\lambda}}^\d(\sqrt{K})=(\sqrt{\lambda}-i\sqrt{|\nu_0(0)|})(\sqrt{\lambda}+i\sqrt{|\nu_0(0)|})f_1(\sqrt{\lambda})$$
  for some entire function $f_1$.
Thus,
\begin{align}
c^\d(\lambda,0,\pi)&= \pi K
J_{\sqrt{\lambda}}^\d(\sqrt{K})J_{-\sqrt{\lambda}}^\d(\sqrt{K})\nonumber\\
&= \pi
K(\lambda-\nu_0(0))^2f_1(\sqrt{\lambda})f_1(-\sqrt{\lambda}),\nonumber
\end{align}
implying $\nu_1(0)<0$. But then (\ref{trace_for1}) would diverge, since
there exists at least one $\nu_k(0)$ in the  open interval
$((n-1)^2,n^2)$ for every $n\in\N$. This is a contradiction, and hence
there is no negative Neumann eigenvalue.

The previous argument also shows that if
$c^\d(n_0^2,0,\pi)=0$ for some $n_0\in \N$, then $c^\d(\lambda,0,\pi)$
has at least a double  zero at $\lambda=n_0^2$, while
$c^\d(\lambda,0,\pi)$ could have a simple zero at $\lambda=0$.

Moreover, still assuming $0<K\leq 1$, if there were more than one
Neumann eigenvalue $\nu_k(0)$ in some interval $((n-1)^2,n^2)$, there
would be at least three $\nu_k(0)$ in this interval, counting
multiplicity since $\nu_k(0)\not=m^2$. This would violate the convergence
of the sum in (\ref{trace_for1}). Hence, there is exactly one $\nu_k(0)$
in each interval $((k-1)^2,k^2)$, $k\in\N$. Since $\nu_k(0)-k^2>0$
for all $k\in\N_0$, we infer that
\begin{align}
&0<\nu_0(0)<\frac{K}{2}, \nonumber \\
&k^2<\nu_k(0)<k^2+\frac{K}{2}-\nu_0(0)
-\sum_{m=1}^{k-1}\nu_m(0)<k^2+\frac{K}{2}-\nu_0(0), \nonumber \\
& \hspace*{9.8cm} k\in\N.  \nonumber
\end{align}

\noindent{\it Proof of (i) and (iii).}  Suppose that $K>1$. We recall that
$\nu_k(x_0)$ is numbered according to $k\in\N_0$. In proving Theorem
\ref{mu_real} (i) and (iii), we used Lemma \ref{lemma_bessel} (6), while
here we need Lemma \ref{lemma_bessel} (7). The point of Lemma
\ref{lemma_bessel} (7) in the proof is that if $c^\d(n^2,0,\pi)=0$ then
$c^\d(k^2,0,\pi)\not=0$ for $k=(n-1)$ and $k=(n+1)$.  The rest of the
proof is similar to that of Theorem \ref{mu_real} (i) and (iii).
Hence, we omit further details.
\end{proof}

Next, we investigate $\mu_j(\pi/2)$, $j\in\N$.
\begin{theorem}\label{pi_2}
Suppose that $V(x)=K e^{2ix}$ and $K>0$.
Then,
$$(2j-1)^2<\mu_{2j-1}(\pi/2)\leq
\mu_{2j}(\pi/2)<(2j)^2, \quad j\in\N,$$
and hence $\mu_j(\pi/2)\not=E_m$, $m\in\N_0$.

Moreover, if $0<K\leq 1$, then
$$(2j-1)^2<\mu_{2j-1}(\pi/2)<
(2j-\frac{1}{2})^2<\mu_{2j}(\pi/2)<(2j)^2, \quad
j\in\N.$$
\end{theorem}
\begin{proof} By (\ref{s_eq}) and (\ref{trans_eq}), $\mu_j(\pi/2)$ are
the values of $\lambda$ for which
$$s(\lambda,\pi/2,3\pi/2)=\pi
J_{-\sqrt{\lambda}}(e^{i\pi/2}\sqrt{K})J_{\sqrt{\lambda}}(e^{i\pi/2}\sqrt{K})=0.$$
If $\lambda\geq 0$, then all the zeros $u$ of $J_{\sqrt{\lambda}}(u)$ are
real and one concludes that
$J_{\sqrt{\lambda}}(e^{i\pi/2}\sqrt{K})\not=0$. Thus, the zeros of
$s(\lambda,\pi/2,3\pi/2)$ agree with those of
$J_{-\sqrt{\lambda}}(i\sqrt{K})$ and hence with those of
$\left[\left(i\sqrt{K}\,
2^{-1}\right)^{\sqrt{\lambda}}J_{-\sqrt{\lambda}}(i\sqrt{K})\right]$.

Let
\begin{equation}\label{f_est}
f(\alpha)=\left(\frac{i\sqrt{K}}{2}\right)^{-\alpha}J_{\alpha}(i\sqrt{K})=\sum_{m=0}^{\infty}\frac{\left(\frac{\sqrt{K}}{2}\right)^{2m}}{m!\Gamma(m+\alpha+1)}.
\end{equation}
We will show below that if $0<K\leq 1$, then
\begin{equation}\label{f_ineq}
f(-n)>0 \text{ and }  f(-(2n-1/2))<0 \text{ for all $n\in\N.$}
\end{equation}
By the intermediate value theorem this then implies the existence of at
least one $\mu_j(\pi/2)$ in every interval of the form
$$(-(2n-1)^2, -(2n-1/2)^2),\quad  (-(2n-1/2)^2, -(2n)^2), \quad n\in\N.$$
Thus, the convergence of the trace formula
\begin{equation}
-K=Ke^{2ix_0}\Big|_{x_0=\pi/2}=E_0
+\sum_{j=1}^{\infty}\big(E_{2j-1}+E_{2j}-2\mu_j(\pi/2)\big)
\end{equation}
implies that there is exactly one eigenvalue $\mu_j(\pi/2)$ in each
interval above, and these are all the $\mu_j(\pi/2)$.

Since $J_{-n}(i\sqrt{K})=(-1)^n J_n(i\sqrt{K})$, $n\in\Z$, and
since by (\ref{f_est}), $f(n)>0$, $n\in\N_0$, we infer from Lemma
\ref{lemma_bessel} (4) that
\begin{align}
f(-n)&=(i\sqrt{K}/2)^nJ_{-n}(i\sqrt{K})\nonumber \\
&=(-1)^n(i\sqrt{K}/2)^{2n}f(n)>0,\quad n\in\N_0. \label{fn_est}
\end{align}

Next, we show that $f(-2n+1/2)<0$, $n\in\N$. For every $0\leq m\leq 2n-2$,
repeated use of the formula
$\Gamma(\lambda)=\frac{1}{\lambda}\Gamma(\lambda+1)$, implies
\begin{equation}\label{gamma_eq}
\Gamma(m-2n+3/2)
=\frac{(-2)^{2n-m-1}\Gamma(1/2)}{\prod_{j=1}^{2n-m-1}(4n-2m-(2j+1))}.
\end{equation}
  Thus, from (\ref{f_est}),
\begin{align}
f(-2n+1/2)
&=\sum_{m=0}^{\infty}\frac{\left(\frac{\sqrt{K}}{2}
\right)^{2m}}{m!\Gamma(m-2n+1/2+1)}\nonumber\\
&=\sum_{m=0}^{2n-2}(-1)^{m+1}
\frac{\prod_{j=1}^{2n-m-1}(4n-2m-(2j+1))}{2^{2n-1}m!
\Gamma(1/2)}\left(\frac{K}{2}\right)^{m}\nonumber\\
&\quad+\sum_{m=2n-1}^{\infty}\frac{\left(\frac{\sqrt{K}}{2}
\right)^{2m}}{m!\Gamma(m-2n+1/2+1)}. \label{sec_last}
\end{align}

Moreover,
\begin{equation}
\sum_{m=2n-1}^{\infty}\frac{\left(\frac{\sqrt{K}}{2}\right)^{2m}}{m!\Gamma(m-2n+1/2+1)}
<\frac{\left(\frac{K}{4}\right)^{2n-1}e^{K/4}}{(2n-1)!\sqrt{\pi}}, \quad
n\in\N,
\end{equation}
and hence, if $0<K\leq 1$,
\begin{align}
f(-2n+1/2)
&< -\frac{\left(\frac{K}{2}\right)^{2n-2}}{
2^{2n-1}(2n-2)!\sqrt{\pi}}+\frac{\left(\frac{K}{4}\right)^{2n-1}e^{K/4}}{(2n-1)!\sqrt{\pi}}\nonumber\\
&=-\frac{\left(\frac{K}{2}\right)^{2n-2}}{2^{2n-1}\sqrt{\pi}(2n-2)!}\left[1-
\frac{\left(\frac{K}{2}\right)e^{K/4}}{(2n-1)}\right]<0. \label{f_half}
\end{align}
Here  we used the fact that the sign of an alternating sum with terms of
decreasing magnitudes agrees with the sign of the first term.

In order to prove the first part of the theorem, we note that
(\ref{fn_est}) holds for every $K>0$, while (\ref{f_half}) may not hold
for some $K>1$. Next, for each $K>1$ and $\frac{1}{K}\leq \varepsilon\leq
1$ we introduce the family of potentials
$V(x)=\varepsilon K e^{2ix}$. Again, we use the notation
$\mu_j(\varepsilon K,0)$  to indicate the $\varepsilon$-dependence of these
eigenvalues. Then, for each $j\in\N$, the function $\varepsilon\mapsto
\mu_j(\varepsilon K,0)$ is  continuous since $J_{\nu}(u)$ is an entire
function of $\nu$ for each fixed $u\not=0$, and an analytic function of
$u$ on the positive real axis.

When $\varepsilon=\frac{1}{K}$, that is, when $V(x)=e^{2ix}$, we proved above
that
$$(2j-1)^2<\mu_{2j-1}(1,\pi/2)\leq \mu_{2j}(1,\pi/2)<(2j)^2,
\quad j\in\N.$$
So as $\varepsilon$ increases to $1$, the real numbers $\mu_j(\varepsilon
K,\pi/2)$ cannot become $(2j-1)^2$ or $(2j)^2$ due to
(\ref{fn_est}) which holds for all $K>0$. This completes the proof.
\end{proof}

Next, we investigate $\nu_k(\pi/2)$ for $k\in\N_0$.
\begin{theorem}
Suppose that $V(x)=K e^{2ix}$ and $K>0$.
Then
$$(2k)^2<\nu_{2k}(\pi/2)\leq \nu_{2k+1}(\pi/2)<(2k+1)^2, \quad
k\in\N_0,$$
and hence $\nu_k(\pi/2)\not=E_m$, $k\in\N_0$.

Moreover, if $0<K\leq \frac{1}{2}$, then
$$(2k)^2<\nu_{2k}(\pi/2)<
(2k+\frac{1}{2})^2<\mu_{2k+1}(\pi/2)<(2k+1)^2, \quad k\in\N_0.$$
\end{theorem}
\begin{proof}
One can follow the arguments in the proof of Theorem \ref{pi_2}.
By (\ref{c_eq}) and (\ref{trans_eq}),   $\nu_k(\pi/2)$ are values of
$\lambda$ for which
$$c^\d(\lambda,\pi/2,3\pi/2)=i\pi
J_{-\sqrt{\lambda}}^\d(i\sqrt{K})J_{\sqrt{\lambda}}^\d(i\sqrt{K})=0.$$
If $\lambda\geq 0$, then $J_{\sqrt{\lambda}}^\d(i\sqrt{K})\not=0$ and
hence the zeros of $c^\d(\lambda,\pi/2,3\pi/2)$ agree with those of
$J_{-\sqrt{\lambda}}^\d(i\sqrt{K})$ and hence they agree with those of
$\left(\frac{i\sqrt{K}}{2}\right)^{\sqrt{\lambda}+1}
J_{-\sqrt{\lambda}}^\d(i\sqrt{K})$. Next, define
\begin{equation}\nonumber
g(\alpha)=\left(\frac{i\sqrt{K}}{2}\right)^{-\alpha+1}J_{\alpha}^\d(i\sqrt{K}),\quad
\alpha\in\R.
\end{equation}
  Then one can show that for each $K>0$,
\begin{equation}
g(n)>0,\quad n\in\Z.
\end{equation}
Here the identity
$$J_{\nu}^\d(u)=\frac{\nu}{u}J_{\nu}(u)-J_{\nu+1}(u)$$
(see \cite[pp.\  45]{WATSON}) turns out to be useful.

Moreover, one can show that if $0<K\leq \frac{1}{2}$, then
\begin{equation}\nonumber
g(-2n-1/2)<0,\quad n\in\N_0.
\end{equation}
  Thus, if $0<K\leq \frac{1}{2}$, then $g(\alpha)$ has at least one zero in
each interval of the form $(-2n-1,-2n-1/2)$ and $(-2n-1/2,-2n)$,
$n\in\N_0$. Using the trace formula (\ref{Neumann_trace1}) one can
show that there are no other zeros and hence this proves the second part
of the theorem.

The proof of the first part of the theorem is analogous to that of
Theorem \ref{pi_2}.
\end{proof}

\begin{theorem}
Suppose that $V(x)=K e^{2ix}$ and $K>0$. Then
for every $x\in (0,\pi)$, $\mu_j(x),\,\nu_k(x)\not=E_m$.
\end{theorem}
\begin{proof}
{} From (\ref{trans_eq}) one knows that $\mu_j(x)$ are zeros
of
\begin{equation}\label{J_eq3}
s(\lambda,x,x+\pi)=\pi
J_{\sqrt{\lambda}}(e^{ix}\sqrt{K})J_{-\sqrt{\lambda}}(e^{ix}\sqrt{K}),
\end{equation}
and that $\nu_k(x)$ are zeros of
\begin{equation}\label{J_eq4}
c^\d(\lambda,x,x+\pi)=\pi\sqrt{K}e^{ix}
J_{\sqrt{\lambda}}^\d(e^{ix}\sqrt{K})
J_{-\sqrt{\lambda}}^\d(e^{ix}\sqrt{K}).
\end{equation}
Since all zeros of $J_n(u)$ and  $J_n^\d(u)$
are real for $n\geq 0$, we see that if $0<x<\pi$,
then $\mu_j(x),\,\nu_k(x)\not=n^2$, $n\in\N_0$, because $e^{ix}$ is
non-real.
\end{proof}

\noindent {\bf Acknowledgments.} The author thanks Fritz Gesztesy for
initiating this research, for many pertinent references on this subject,
and for showing him Percy Deift's unpublished manuscript. He also
thanks Mark Ashbaugh for some help with references on Bessel functions.

{\sc e-mail:}  kcshin@math.missouri.edu
\end{document}